\documentclass[12pt]{iopart}

\newcommand{\pdf}{\textsl{pdf}}
\newcommand{\chf}{\textsl{chf}}
\newcommand{\id}{\textsl{id}}
\newcommand{\SDE}{\textsl{SDE}}
\newcommand{\PIDE}{\textsl{PIDE}}
\newcommand{\GH}{\mathcal{GH}}
\newcommand{\VG}{\mathcal{VG}}
\newcommand{\ST}{\mathcal{T}}
\newcommand{\rv}{\textsl{rv}}
\newcommand{\binom}[2]{{{#1}\choose {#2}}}
\newcommand{\tod}{\stackrel{d}{\longrightarrow}}
\newcommand{\ci}{\mathrm{c i}\,}
\newcommand{\si}{\mathrm{s i}\,}

\newtheorem{prop}{Proposition}[section]

\usepackage{iopams}
\usepackage{setstack}
\usepackage{graphicx}

\begin{document}

\title[N Cufaro Petroni: Mixtures in non stable L\'evy processes]{Mixtures in non stable L\'evy processes}

\author{N Cufaro Petroni}

\address{Department of Mathematics and TIRES, Universit\`a di Bari;
INFN Sezione di Bari, \\
via E Orabona 4, 70125 Bari, Italy}
\ead{cufaro@ba.infn.it}

\begin{abstract}

We analyze the L\'evy processes produced by means of two
interconnected classes of non stable, infinitely divisible
distribution: the Variance Gamma and the Student laws. While the
Variance Gamma family is closed under convolution, the Student one
is not: this makes its time evolution more complicated. We prove
that -- at least for one particular type of Student processes
suggested by recent empirical results, and for integral times -- the
distribution of the process is a mixture of other types of Student
distributions, randomized by means of a new probability
distribution. The mixture is such that along the time the asymptotic
behavior of the probability density functions always coincide with
that of the generating Student law. We put forward the conjecture
that this can be a general feature of the Student processes. We
finally analyze the Ornstein--Uhlenbeck process driven by our L\'evy
noises and show a few simulation of it.

\end{abstract}




\section{Introduction}\label{intro}

Since a few years the L\'evy processes enjoy considerable popularity
in several different fields of research from statistical physics to
mathematical finance (Paul and Baschnagel 1999, Mantegna and Stanley
2001, Barndorff--Nielsen \etal 2001 and Cont and Tankov 2004 are
just a few examples of books reviewing the large body of literature
on this subject). In the former field, however, the interest has
been generally confined to $\alpha$--stable processes which are an
important particular sub--class of L\'evy processes (Bouchaud and
Georges 1990, Metzler and Klafter 2000, Paul and Baschnagel 1999,
Woyczy\'nski 2001), while studies about non stable, infinitely
divisible L\'evy processes abound mainly in the latter field (see
for example Cont and Tankov 2004 and references quoted therein). The
appeal of the $\alpha$--stable distributions is justified by the
properties of scaling and self--similarity displayed by the
corresponding processes, but it must also be remarked that these
distributions show a few features that partly impair their
usefulness as empirical models. First of all the non gaussian stable
laws always have infinite variance. This makes them rather suspect
as a realistic tool and prompts the introduction of
\textsl{truncated} stable distributions which, however, are no
longer stable. Then the range of the $x$ decay rates of the
probability density functions can not exceed $x^{-3}$, and this too
introduces a particular rigidity in these models. On the other hand
the more general L\'evy processes are generated by infinitely
divisible laws and do not necessarily have these problems, but they
can be more difficult to analyze and to simulate. Beside the fact
that they do not have natural scaling properties, the laws of their
increments could be explicitly known only at one time scale. In fact
their time evolution is always given in terms of characteristic
functions, but the marginal densities may not be calculable. This is
a feature, however, that they share with most stable processes,
since the probability density functions of the non gaussian stable
laws are explicitly known only in precious few cases.

The need to go beyond the processes generated by stable
distributions stems also from other recent advances in the field of
the fractional differential equations. The evolution equations of
the L\'evy processes can be put in terms of pseudo--differential
operators whose symbols are just the characteristic exponents of the
processes (Jacob and Schilling 2001, Cont and Tankov 2004). The most
popular form taken by these equations is that of the fractional
differential equations, and this generalization of the diffusion
equations can be put in connection with L\'evy noises with
non--Gaussian stable distributions (Gorenflo and Mainardi 1998a,
1998b, Metzler and Klafter 2000). It has been put in evidence in a
few papers (Chechkin \etal 2003, 2004), however, that in the case of
L\'evy flights confined by symmetric quartic potentials the
stationary probability density functions show two unexpected
properties: in fact not only they are bimodal, but they also have a
finite variance, differently from what happens to the non--Gaussian,
stable law of the system noise. This suggests that, under particular
dynamical conditions, the stochastic evolutions produced by stable
L\'evy noises end up in non stable distributions, and hence hints to
a new physical interest beyond the pale of the stable laws.

Some new applications for the L\'evy, infinitely divisible but not
stable processes begin also to emerge in other physical domains
(Cufaro Petroni \etal 2005, 2006, Vivoli \etal 2006): as we will see
in the following the statistical characteristics of some recent
model of the collective motion in the charged particle accelerator
beams seem to point exactly in the direction of some kind of Student
infinitely divisible process. At the present stage of our inquiry
the proposed model for the particle beams is only phenomenological
and it lacks a complete, underlying, physical mechanism producing
the noise. This however brings to the fore the problem of the
dynamical description of complex systems. The infinitely divisible
L\'evy processes with a jump component are indeed interesting also
in the light of the connection established between Markov processes
and quantum phenomena by the stochastic mechanics. This latter is a
model universally known for its original application to the problem
of building a classical stochastic model for quantum mechanics
(Nelson 1967, 1985, Guerra 1981, Morato 1982, Guerra and Morato
1983), but in fact it is a very general model which is suitable for
a large number of stochastic dynamical systems (Albeverio, Blanchard
and H\o gh-Krohn 1983, Paul and Baschnagel 1999, Cufaro Petroni
\etal 1999, 2000, 2003, 2004). As recently proposed, a stochastic
mechanics with jumps driven by a non gaussian L\'evy process could
find applications in the physical and technological domain (Cufaro
Petroni \etal 2005, 2006). The presence of jumps could for instance
be instrumental in building reasonable models for the formation of
halos in beams of charged particles in accelerators. On the other
hand this would not be the first time that L\'evy processes find
applications in quantum theory since they have already been used to
build models for spinning particles (De Angelis and Jona--Lasinio
1982), for relativistic quantum mechanics (De Angelis 1990), and in
stochastic quantization (Albeverio, R\"udiger and Wu 2001).

The standard way to build a stochastic dynamical system is to modify
the phase space dynamics by adding a Wiener noise $\mathbf{B}(t)$ to
the momentum equation only, so that the usual relations between
position and velocity are preserved:
\begin{equation*}
m\,d{\bf Q}(t)={\bf P}(t)\,dt\,,\qquad\quad d{\bf P}(t)={\bf
F}(t)\,dt+\beta\, d{\bf B}(t)\,.
\end{equation*}
In this way we get a derivable, but non Markovian position process
${\bf Q}(t)$. An example of this approach is that of a Brownian
motion in a fluid described by an Ornstein--Uhlenbeck system of
stochastic differential equations. Alternatively we can add a Wiener
noise ${\bf W}(t)$ with diffusion coefficient $D$ directly to the
position equation:
\[
d{\bf Q}(t)\;=\;{\bf v}_{(+)}({\bf Q}(t),t)\,dt+\sqrt{D}\,d{\bf
W}(t)\,.
\]
and get a Markovian, but not derivable ${\bf Q}(t)$. In this way the
stochastic system is reduced to a single stochastic differential
equation since we are obliged to drop the second (momentum)
equation. The standard example of this reduction is the Smoluchowski
approximation of the Ornstein--Uhlenbeck process in the overdamped
case. As a consequence we will now work only in a configuration, and
not in a phase space; but this does not prevent us from introducing
a dynamics either by generalizing the Newton equations (Nelson 1967,
1985, Guerra 1981), or by means of a stochastic variational
principle (Guerra and Morato 1983). From this stochastic dynamics,
which now notably enjoys a measure of time--reversal invariance, two
coupled equations can be derived which are equivalent to a
Schr\"odinger equation, prompting the idea of a stochastic
foundation of quantum mechanics. In fact the stochastic mechanics
can be used to describe more general stochastic dynamical systems
satisfying fairly general conditions: it is known since longtime
(Morato 1982), for example, that for any given diffusion there is a
correspondence between diffusion processes and solutions of this
kind of Schr\"odinger equations where the Hamiltonians come from
suitable vector potentials. The usual Schr\"odinger equation, and
hence true quantum mechanics, is recovered when the diffusion
coefficient coincides with $\hbar/2m$, namely is connected to the
Planck constant. However we are interested here not only in a
stochastic model of quantum mechanics, but also to the general
description of complex systems as a particle beams, and to this end
it would be very interesting -- as already remarked -- to be able to
generalize the stochastic mechanical scheme to the case of non
Gaussian L\'evy noises. The road to this end, however, is fraught
with technical difficulties, so that a better understanding of the
possible underlying L\'evy noises should be considered as a first,
unavoidable step.

In this light the aim of this paper is to study a few examples of
non stable, infinitely divisible processes, and in particular we
will focus our attention on the Student processes. Since the Student
family of laws is infinitely divisible but non closed under
convolution the process distribution will not be Student at every
time. We will show however that, at least in particular cases, the
process transition law is a mixture of a finite number of Student
laws, and it is suggested that this could be a general feature of
the Student processes. On the other hand it can also be seen that
for every finite time the spatial asymptotic behavior always is the
same as that of the Student distribution at the characteristic time
scale; and this turns out to be exactly the behavior put in evidence
by Vivoli \etal 2006 in the solutions of the complex dynamical
system used to study the behavior of beams of charged particles in
accelerators.

We will limit our considerations to one dimensional models without
going into the problem of the dependence structure of a multivariate
process (see for example Cont and Tankov 2004), and we will not
pretend any completeness or generality: our aim is rather to present
the features of a few selected processes to gain a deeper insight
into their possible general behavior. The paper is organized as
follows: in Section~\ref{lproc} we recall a few, well known facts
about the L\'evy processes and in particular the connection between
the transition function $p(x,t|\,y,s)$ and a triplet of functions
$A(y,s)$, $B(y,s)$ and $W(x|\,y,s)$ characteristic of a L\'evy
process. We also propose a different simplified, heuristic procedure
to find the explicit form of $A$, $B$ and $W$: a procedure not
completely general, but which works well enough for the rather
regular transition functions discussed in this paper. In the
Section~\ref{GHlaws} we analyze the behavior of two families of laws
(the Variance Gamma, and the Student laws) which are particular
limit cases of a larger class of infinitely divisible laws: that of
the Generalized Hyperbolic laws which received considerable
attention in recent years (Raible 2000, Eberlein and Raible 2000,
Eberlein 2001, Cont and Tankov 2004, and references quoted therein).
Our two families are in a certain sense conjugate to each other
since the roles of their probability density functions and
characteristic functions are interchanged. Let us remark here that
all the laws that we take in consideration in this paper are
infinitely divisible, but -- with a few notable exceptions -- not
stable. We then pass in Section~\ref{VGproc} to study the L\'evy
process produced by the Variance Gamma distributions: since this
class is closed under convolution, it will be easy enough to find
both the characteristic triplet, and the laws of the increments for
every value of the time interval. Apparent similarities
notwithstanding, the case of the Student processes discussed in
Section~\ref{STproc} is rather different from the previous one. In
fact the Student family is not even closed under convolution so that
we do not have explicit expressions for the transition laws at every
time scale. As a consequence we will restrict our attention to a
subclass of Student processes by choosing particular (but not
trivial) values for the parameters, and we will get results about
(a) the spatial asymptotic behavior of the transition functions at
every time, (b) the explicit form of the transition functions at
time intervals which are integral multiples of a characteristic time
constant, and (c) the form of the L\'evy triplet of functions. In
particular we will find that, at discrete times, the process law
turns out to be a mixture of a finite number of Student
distributions by means of a new kind of time dependent discrete
probability distribution. We finally discuss in
Section~\ref{pathprop} some pathwise properties of our non stable
processes by showing also a few simulations of the
Ornstein--Uhlenbeck processes driven by our L\'evy noises, and
conclude with some remark about the perspectives of future research.

\section{L\'evy processes generated by \id\ laws}\label{lproc}

A L\'evy process $X(t)$ is a stationary, stochastically continuous,
independent increment Markov process. It is well known that the
simplest way to produce its transition laws is to start with a
\textsl{type} of infinitely divisible (\id)\ distributions (see
Gnedenko and Kolmogorov 1968, Lo\`eve 1978, 1987 and Sato 1999 for a
more recent monograph): if we focus our attention on centered laws,
a type of these \textsl{generating} laws can be given by the family
of their characteristic functions (\chf) $\varphi(au)$ with a
spatial scale parameter $a>0$. The \chf\ of the transition law of
our stationary process in the time interval $[s,t]$ will then be
\begin{equation}\label{procchf}
    \Phi(au,t-s)=\left[\varphi(au)\right]^{(t-s)/T}
\end{equation}
where $T$ is a suitable constant playing the role of a time scale
parameter, while the transition probability density function (\pdf)\
with initial condition $X(s)=y,\;\mathbb{P}$-q.o.\ will be recovered
by an inverse Fourier transform
\begin{eqnarray}
  p(x,t|\,y,s) &=& \frac{1}{2\pi}\lim_{M\to+\infty}\int_{-M}^M\Phi(au,t-s)\,e^{-i(x-y)u}\,du \nonumber\\
   &=&\frac{1}{2\pi}\lim_{M\to+\infty}\int_{-M}^M\left[\varphi(au)\right]^{(t-s)/T}\,e^{-i(x-y)u}\,du\label{invft}
\end{eqnarray}
and -- because of stationarity -- will only depend on the
differences $x-y$ and $t-s$.

The parameters $a$ and $T$ play a role in the \textsl{scale
invariance} properties of the process. When the generating family of
\id\ laws is closed under convolution the transition laws remain
within this same family all along the evolution, and the changes are
summarized just in a time dependence of some parameter of the \pdf.
But in the case of \textsl{stable} laws there is more. If for
instance -- as in the Wiener process -- the generating type of law
is the normal, centered $\mathcal{N}(0,a)$ it is well known that the
transition law (with $y=0$ and $s=0$ for simplicity) is just
$\mathcal{N}(0,a\sqrt{t/T})$, namely it is always normal, but with a
time dependent parameter: the variance, changing linearly with the
time as $Dt$, where $D=a^2/T$ is the diffusion coefficient. This
means that the overall behavior of the process is ruled only by $D$,
and not by $a$ and $T$ separately. As a consequence the particular
values of $a$ and $T$, namely the particular units of measurement,
are immaterial and we have the scale invariance. This gives to the
Wiener process its property of \textsl{self--similarity}: no matter
at what space--time scale (namely irrespectively to the values you
give to $a$ and $T$, provided that $D=a^2/T$ keeps the same value)
you choose to observe the process, the trajectories always will look
the same.

These properties of the Wiener process are shared by all the other
L\'evy processes generated by stable -- even non normal -- laws, but
not in general by the processes generated by other, non stable \id\
laws. It must be remarked, however, that all the non gaussian stable
laws do not have a finite variance, and show a rather restricted
range of possible decays for large $x$: features that partly impair
a realistic use of them in empirical situations. On the other hand
families of non stable, \id\ laws can still be closed under
convolution, as it is for instance the case of the compound Poisson
laws $\mathcal{P}(\lambda,a;\chi)$ with \chf\
$\varphi(au)=e^{\lambda[\chi(au)-1]}$, where $\chi(u)$ is the \chf\
of the jump distribution. This means again that the evolution of the
transition law of a compound Poisson process can always be
summarized in the time dependence of the Poisson parameter as
$\mathcal{P}(\lambda t/T,a;\chi)$, but with respect to the Wiener
case there are important differences: while all the transition laws
of a Wiener process belong to the same (normal) type, Poisson
transition laws with different parameters do not. The normal laws
are indeed stable, while the Poisson laws are only \id, and Poisson
laws with different values of $\lambda$ do not belong to the same
type. Moreover, while a change in the $T$ value can always be
compensated by a corresponding change of $\lambda$ so that
$\lambda/T$ remains the same, the roles of $a$ and $T$ in a compound
Poisson process, at variance with the Wiener case, remain completely
separated and we do not have the same kind of self--similarity.

The less simple case of processes is finally that generated by
families of \id\ laws which are not even closed under convolution,
since in this event the transition distributions do not remain
within the same family, and the overall evolution can not be
summarized just in the time dependence of some parameter. As we will
see in the following this is far to be an uncommon situation and
this paper is mainly devoted to the analysis of particular processes
of this kind. It must be kept in mind that in this last case the
role of the scale parameters becomes relevant since a change in
their values can no longer be compensated by reciprocal changes in
other parameters. This means that, to a certain extent, a change in
these scale constants produces different processes, so that for
instance we are no longer free to look at the process at different
time scales by presuming to see the same features. We should remark,
on the other hand, that -- at variance with the stable, non gaussian
case -- the \pdf's of the \id\ distributions can have both a wide
range of decay laws for $|x|\to+\infty$, and a finite variance
$\sigma^2$. For these L\'evy processes generated by \id\ laws with
finite variance $\sigma^2$ it is finally easy to see that -- due to
the fact that the process has independent increments -- the variance
always is finite and grows linearly with the time as $\sigma^2t/T$:
a feature typical of the ordinary (non anomalous) diffusions.

\subsection{The decomposition of a L\'evy process}

The evolution equations of a process driven by a L\'evy noise can be
given either as partial integro--differential equations (\PIDE)\ for
the transition functions of the process (Lo\`eve 1978, Gardiner
1997), or as stochastic differential equations (\SDE)\ for its
trajectories (Applebaum 2004, {\O}ksendal and Sulem 2005, Protter
2005). In both cases the structure of the evolution is given in
terms of some characteristic triplet of functions. For simple L\'evy
process of course this triplet will give rise just to its L\'evy
decomposition in a drift, a Brownian and a jump term. In this paper
we will choose to follow the description in terms of \PIDE, and it
is important to recall how this characteristic triplet is related to
the transition functions. We will not attempt to give here a
complete and rigorous survey of the argument, but we will limit
ourselves to fix the notation in a rather simplified form (see for
example Gardiner 1997, but also for a more rigorous approach
L\'eandre 1987, Ishikawa 1994, Sato 1999, Barndorff--Nielsen 2000
and R\"uschendorf and Woerner 2002) suitable for the cases that we
will analyze. In particular we suppose to consider only processes
endowed with well behaved \pdf's, so that (apart from an initial
distribution) the process is completely defined by its transition
\pdf\ $p(x,t|\,y,s)$. If then we define the triplet of functions
\begin{eqnarray}
  A(y,s) &=& \lim_{\epsilon\to0^+}\lim_{\Delta t\to0}\frac{1}{\Delta t}
                            \int_{|x-y|<\epsilon}(x-y)p(x,s+\Delta t|\,y,s)\,dx\label{Afunct} \\
  B(y,s) &=& \lim_{\epsilon\to0^+}\lim_{\Delta t\to0}\frac{1}{\Delta t}
                            \int_{|x-y|<\epsilon}(x-y)^2p(x,s+\Delta t|\,y,s)\,dx\label{Bfunct} \\
  W(x|\,y,s) &=& \lim_{\Delta t\to0}\frac{p(x,s+\Delta t|\,y,s)}{\Delta
  t}\,,\qquad\qquad x\neq y\label{Wfunct}
\end{eqnarray}
it can be seen that the \pdf's of the process satisfy the following
(\textsl{forward}) \PIDE
\begin{eqnarray}\label{pide}
    \partial_tp(x,t)&=&-\partial_x[A(x,t)p(x,t)]+\frac{1}{2}\,\partial^2_x[B(x,t)p(x,t)]\nonumber\\
             & & \quad+\lim_{\epsilon\to0^+}\int_{|x-z|\geq\epsilon}
                \left[W(x|\,z,t)p(z,t)-W(z|\,x,t)p(x,t)\right]\,dz
\end{eqnarray}
the transition \pdf\ being the solution corresponding to the initial
condition $p(x,s^+|\,y,s)=\delta(x-y)$. In the case of stationary
processes (as our L\'evy processes are) the transition \pdf\
$p(x,t|\,y,s)$ depends on its variables only trough their
differences $x-y$ and $t-s$. As a consequence $A$ and $B$ are simply
constants, while $W(x|\,y,s)=W(x-y)$. It is also known that $A$
plays the role of a drift coefficient, while $B$ is a diffusion
coefficient connected to the Brownian component of the process;
finally $W(x|\,y,s)$, defined only for $x\neq y$, is the density of
the L\'evy measure of the process. The knowledge of the
characteristic triplet is also instrumental to write down the \PIDE\
(or alternatively the \SDE)\ for other processes driven by a L\'evy
noise.

In order to calculate the characteristic triplet of a L\'evy process
decomposition from~\eref{Afunct},~\eref{Bfunct} and~\eref{Wfunct} we
are supposed to explicitly know its transition \pdf. We will see in
the following, however, that given the \chf's of an \id\
distribution it is very easy to write the \chf~\eref{procchf} of the
process increments, but also that in general it is not a simple task
to explicitly calculate the transition \pdf\ by the inverse Fourier
transform~\eref{invft}. We then propose here a different procedure
to calculate $A, B$ and $W$ directly from the process \chf\ which is
surely a known quantity for a L\'evy process, by adding however that
at the present stage its derivation is only heuristic. To this end
let us remark that from~\eref{procchf} and~\eref{invft} the
transition \pdf\ will have the form
\begin{equation*}
    p(x,s+\Delta t|\,y,s)=\frac{1}{2\pi}\lim_{M\to+\infty}\int_{-M}^M
             \left[\varphi(au)\right]^{\Delta t/T}e^{-iu(x-y)}du
\end{equation*}
so that, by supposing (which is fair for all the cases that we will
consider in this paper) $\varphi(-\infty)=\varphi(+\infty)=0$, we
get with an integration by parts
\begin{eqnarray*}
  \lefteqn{\frac{p(x,s+\Delta t|\,y,s)}{\Delta t}}\\
      && \qquad\qquad = \frac{a}{2\pi
             i(x-y)T}\lim_{M\to+\infty}\int_{-M}^M
             \left[\varphi(au)\right]^{\Delta
             t/T}\frac{\varphi\,'(au)}{\varphi(au)}\,e^{-iu(x-y)}du
\end{eqnarray*}
If now we suppose that our functions are regular enough to allow
both to exchange the two limits for $\Delta t\to0$ and for
$M\to+\infty$, and to perform the limit for $\Delta t\to0$ under the
integral, we immediately have
\begin{equation*}
  W(x|\,y,s) =W(x-y)
     =\frac{a}{2\pi i(x-y)T}\lim_{M\to+\infty}\int_{-M}^M\frac{\varphi\,'(au)}{\varphi(au)}\,e^{-iu(x-y)}du
\end{equation*}
namely with $z=x-y$
\begin{equation}\label{Wchf}
  W(z) =\frac{a}{2\pi izT}\lim_{M\to+\infty}\int_{-M}^M\frac{\varphi\,'(au)}{\varphi(au)}\,e^{-iuz}du
\end{equation}
Remark that in~\eref{Wchf} the limit must be understood in the sense
of the distributions, as can be easily checked by applying the
formula to some well known case (either the Wiener, or the Cauchy
process). What is most interesting with respect to the~\Eref{Wfunct}
is that now we can calculate $W(z)$ directly from $\varphi(au)$,
without explicitly knowing the transition \pdf\ $p(x,t|\,y,s)$.

In the same way for $A$ with an integration by parts we have first
of all that
\begin{eqnarray*}
\fl \lefteqn{\frac{1}{\Delta
t}\int_{|x-y|<\epsilon}(x-y)p(x,s+\Delta
    t|\,y,s)\,dx }  \\
    && \qquad\quad = \;\frac{a}{2\pi iT}\int_{|x-y|<\epsilon}\left[\lim_{M\to+\infty}
                    \int_{-M}^M\left[\varphi(au)\right]^{\Delta
                                      t/T}\frac{\varphi\,'(au)}{\varphi(au)}\,e^{-iu(x-y)}du\right]\,dx
\end{eqnarray*}
then, if again it is allowed to freely exchange limits and
integrals, we have
\begin{eqnarray*}
  \lefteqn{\lim_{\Delta t\to0}\frac{1}{\Delta t}\int_{|x-y|<\epsilon}(x-y)p(x,s+\Delta
    t|\,y,s)\,dx }  \\
   && \qquad\qquad = \; \frac{a}{2\pi iT}\int_{|x-y|<\epsilon}\left[\lim_{M\to+\infty}
                    \int_{-M}^M\frac{\varphi\,'(au)}{\varphi(au)}\,e^{-iu(x-y)}du\right]\,dx\\
   && \qquad\qquad = \; \frac{a}{i\pi T}\lim_{M\to+\infty}
                    \int_{-M}^M\frac{\varphi\,'(au)}{\varphi(au)}\,\frac{\sin
                    u\epsilon}{u}\,du
\end{eqnarray*}
and finally
\begin{equation}\label{Achf}
    A(y,s)=A=\frac{a}{i\pi T}\,\lim_{\epsilon\to0^+}\lim_{M\to+\infty}
                    \int_{-M}^M\frac{\varphi\,'(au)}{\varphi(au)}\,\frac{\sin
                    u\epsilon}{u}\,du
\end{equation}
Here it is understood that the two limits (always in the sense of
distributions) and the integration must be performed in the order
indicated since an exchange will produce a trivial -- and wrong --
result. Remark that when $\varphi(au)$ is an even function (as
happens if the process increments are symmetrically distributed
around zero), then $\varphi\,'(au)/\varphi(au)$ is an odd function,
and hence -- since $u^{-1}\sin u\epsilon$ is even -- we immediately
get $A=0$. This is coherent with the fact that, when the increments
are symmetrically distributed, then we do not expect to have a drift
in the process.

As for the coefficient $B$ the usual integration by parts gives
\begin{eqnarray*}
\fl  \lefteqn{\frac{1}{\Delta
t}\int_{|x-y|<\epsilon}(x-y)^2p(x,s+\Delta
    t|\,y,s)\,dx }  \\
    &&  = \;\frac{a}{2\pi iT}\int_{|x-y|<\epsilon}\left[\lim_{M\to+\infty}
                    \int_{-M}^M\left[\varphi(au)\right]^{\Delta
                                      t/T}\frac{\varphi\,'(au)}{\varphi(au)}\,(x-y)e^{-iu(x-y)}du\right]\,dx
\end{eqnarray*}
so that by exchanging limits and integrals we get
\begin{eqnarray*}
  \lefteqn{\lim_{\Delta t\to0}\frac{1}{\Delta t}\int_{|x-y|<\epsilon}(x-y)^2p(x,s+\Delta
    t|\,y,s)\,dx }  \\
   && \qquad\qquad = \; \frac{a}{2\pi iT}\int_{|x-y|<\epsilon}\left[\lim_{M\to+\infty}
                    \int_{-M}^M\frac{\varphi\,'(au)}{\varphi(au)}\,(x-y)e^{-iu(x-y)}du\right]\,dx\\
   && \qquad\qquad = \; \frac{a}{\pi T}\lim_{M\to+\infty}
                    \int_{-M}^M\frac{\varphi\,'(au)}{\varphi(au)}\,\frac{u\epsilon\cos u\epsilon-\sin
                    u\epsilon}{u^2}\,du
\end{eqnarray*}
and finally our coefficient is
\begin{equation}\label{Bchf}
    B(y,s)=B=\frac{a}{\pi T}\,\lim_{\epsilon\to0^+}\lim_{M\to+\infty}
                    \int_{-M}^M\frac{\varphi\,'(au)}{\varphi(au)}\,\frac{u\epsilon\cos u\epsilon-\sin
                    u\epsilon}{u^2}\,du
\end{equation}
Also in this case we see that for our stationary, independent
increment process this coefficient is a constant independent from
the initial coordinates $y$ and $s$.

The formulas~\eref{Wchf}, \eref{Achf} and~\eref{Bchf} can finally be
checked on two well known (stable) cases to give the correct
characteristic triplets: the Wiener process produced by a normal
distribution $\mathcal{N}(0,a)$ with
\begin{equation}\label{normtriple}
    A=0\,,\qquad\quad B=\frac{a^2}{T}\,,\qquad\quad\;\, W(z)=0
\end{equation}
and the Cauchy process produced by a Cauchy distribution
$\mathcal{C}(a)$ with
\begin{equation}\label{cauchytriple}
    A=0\,,\qquad\quad B=0\,,\qquad\qquad W(z)=\frac{a}{\pi Tz^2}
\end{equation}
Remark as in these two stable cases the elements of the triplet do
not depend separately on the two (time and space) scale parameters,
but only on a combination of them so that a change in the time scale
can always be compensated by an exchange in the space scale (and
vice versa): a point giving rise to the scale invariance which in
general is not reproduced in non stable processes, as discussed at
the beginning of this section.

\section{A class of infinitely divisible
distributions}\label{GHlaws}

The increment laws of the L\'evy processes analyzed in this paper
are particular (limiting) cases of a larger class of distributions,
that of the Generalized Hyperbolic (GH) distributions (for their
general properties see for example Raible 2000, Eberlein and Raible
2000, Eberlein 2001, Cont and Tankov 2004, and references quoted
therein). The GH distributions constitute a five--parameter class of
\id, absolutely continuous laws with the following \pdf's (for
$x\in\mathbb{R}$) and \chf's
\begin{eqnarray*}
 \fl f(x+\mu) &=& \frac{e^{\beta x}}{\alpha^{2\lambda-1}\delta^{2\lambda}\sqrt{2\pi}}\,
                    \frac{(\delta\sqrt{\alpha^2-\beta^2})^\lambda}{K_\lambda(\delta\sqrt{\alpha^2-\beta^2})}\,
                    (\alpha\sqrt{\delta^2+x^2})^{\lambda-\frac{1}{2}}K_{\lambda-\frac{1}{2}}(\alpha\sqrt{\delta^2+x^2}) \\
 \fl \quad\;\;\; \varphi(u) &=& e^{i\mu u}\frac{(\delta\sqrt{\alpha^2-\beta^2})^\lambda}{K_\lambda(\delta\sqrt{\alpha^2-\beta^2})}\,
                      \frac{K_\lambda(\delta\sqrt{\alpha^2-(\beta+iu)^2})}{(\delta\sqrt{\alpha^2-(\beta+iu)^2})^\lambda}
\end{eqnarray*}
where $\lambda\in\mathbb{R},\;\alpha>0,\;
\beta\in(-\alpha,\alpha),\;\delta>0,\;\mu\in\mathbb{R}$, and
$K_\nu(z)$ are the modified Bessel functions (Abramowitz and Stegun
1968). Apparently $\alpha$ and $\delta$ play the role of scale
parameters, while $\beta$ is a skewness parameter: the \pdf\ is
symmetric when $\beta=0$. On the other hand $\mu$ is just a
centering parameter: since in this paper our attention will be
focused on the symmetric, centered laws, we will always choose
$\beta=0$ and $\mu=0$ and we will consider the more restricted (but
still large enough) class $\GH(\lambda,\alpha,\delta)$ of the
centered, symmetric GH laws with the following \pdf's and \chf's
\begin{eqnarray}
  f_{GH}(x) &=& \frac{\alpha}{(\delta\alpha)^\lambda K_\lambda(\delta\alpha)\sqrt{2\pi}}\,
                    (\alpha\sqrt{\delta^2+x^2})^{\lambda-\frac{1}{2}}K_{\lambda-\frac{1}{2}}(\alpha\sqrt{\delta^2+x^2})\label{GH1} \\
 \varphi_{GH}(u) &=& \frac{(\delta\alpha)^\lambda}{K_\lambda(\delta\alpha)}\,
                      \frac{K_\lambda(\delta\sqrt{\alpha^2+u^2})}{(\delta\sqrt{\alpha^2+u^2})^\lambda}\label{GH2}
\end{eqnarray}
with $\lambda\in\mathbb{R},\;\alpha>0,$ and $\delta>0$.

The GH class contains many relevant particular cases, also for limit
values of the parameters, and its name comes from the fact that it
contains as sub--class with $\lambda=1$ that of the Hyperbolic
distributions called in this way because the logarithm of their
\pdf\ is a hyperbola. The GH distributions are not always endowed
with finite momenta: this fact depends on the parameter values and
must be explicitly assessed for every particular case. On the other
hand they are all \id, and hence they are good starting points to
build L\'evy processes. In general, however, they are not stable
laws, and in fact they are not even closed under convolution: the
sum of two GH random variables (\rv)\ is not a GH \rv. This means
not only that the corresponding processes will not be self--similar,
but also that often it is not easy to find out what the \pdf\ of the
process looks like even if it is well known at one time. Remark that
the GH class is rich enough to contain also as a limit case the
sub--class of the normal laws $\mathcal{N}(\mu, \sigma)$. Indeed it
can be shown that (in distribution)
\begin{equation*}
    \lim_{\delta\to+\infty}\,\lim_{\lambda\to-\infty}\,\lim_{\alpha\to0^+}
            \GH(\lambda,\alpha,\delta)=\mathcal{N}(0,\sigma)
\end{equation*}
provided that $\delta^2/|\lambda|\to2\,\sigma^2$. In the following
we will study the behavior of the processes produced by two other
particular limit sub--classes that, at variance with the normal
distributions, are not stable besides a few exceptions.

\subsection{The Variance Gamma distributions}\label{ssectionVG}

The Variance Gamma (VG) laws (Madan and Seneta 1987, Madan and
Seneta 1990, Madan and Milne 1991, Madan \etal 1998) are obtained
from $\GH(\lambda,\alpha,\delta)$ in the limit for $\delta\to0^+$.
More precisely, since in general
\begin{equation}\label{Kzero}
    K_\nu(z)=K_{-\nu}(z)\sim\left\{
                  \begin{array}{ll}
                    \frac{1}{2}\,\Gamma(\nu)\,(2/z)\,^{\nu}, &\quad \hbox{for $\nu>0\,$,} \\
                    -\log z, &\quad \hbox{for $\nu=0\,$,} \\
                    \frac{1}{2}\,\Gamma(|\nu|)\,(2/z)\,^{|\nu|}, &\quad \hbox{for $\nu<0\,$,}
                  \end{array}
                \right.
 \qquad    z\to0
\end{equation}
we have for $\lambda>0$ that
\begin{equation*}
    \lim_{\delta\to0^+}(\delta\alpha)^\lambda K_\lambda(\delta\alpha)
    =2^{\lambda-1}\Gamma(\lambda)
\end{equation*}
and hence the \pdf's of the centered, symmetric VG laws -- which
constitute the two parameters family $\VG(\lambda,\alpha)$ -- are
\begin{equation}\label{VG1}
    f_{VG}(x)=\frac{2\alpha}{2^\lambda\Gamma(\lambda)\sqrt{2\pi}}\,
    (\alpha|x|)^{\lambda-\frac{1}{2}}\,K_{\lambda-\frac{1}{2}}(\alpha|x|)\,.
\end{equation}
with $\lambda>0$ and $\alpha>0$. As for the corresponding \chf's of
$\VG(\lambda,\alpha)$ it is readily seen from \eref{GH2} and
\eref{Kzero} that they simply reduce to
\begin{equation}\label{VG2}
    \varphi_{VG}(u)=\left(\frac{\alpha^2}{\alpha^2+u^2}\right)^\lambda.
\end{equation}
It is apparent that $\alpha$ plays the role of a scale parameter,
while $\lambda$ classifies the different types of VG laws. For
$\lambda=1$ the \pdf's and \chf's of the $\VG(1,\alpha)$ laws are
\begin{equation*}
    f(x)=\frac{\alpha}{2}\,e^{-\alpha|x|},\qquad\quad\varphi(u)=\frac{\alpha^2}{\alpha^2+u^2}
\end{equation*}
so that $\VG(1,\alpha)$ is nothing but the class of the Laplace
(double exponential) laws $\mathcal{L}(\alpha)$. The \pdf\
\eref{VG1} has an elementary form only when $\lambda=n+1$ is an
integer with $n=0,1,\ldots$ In fact from
\begin{equation}\label{Ksemint}
    K_{n+\frac{1}{2}}(z)=\sqrt{\frac{\pi}{2z}}\,e^{-z}\sum_{j=0}^n\frac{(n+j)!}{j!(n-j)!}\,\frac{1}{(2z)^j}
\end{equation}
it is easy to see that (with $\ell=n-j$) we have
\begin{equation*}
    f_{VG}(x)=
    \frac{\alpha}{2^{2n+1}}\,e^{-\alpha|x|}\sum_{\ell=0}^n\binom{2n-\ell}{n}\frac{(2\alpha|x|)^\ell}{\ell!}
         \,,\qquad \lambda=n+1=1,2,\ldots
\end{equation*}
For $\lambda\to0$ \Eref{VG2} shows also that our VG laws converge in
law to a distribution degenerate in $x=0$. From the asymptotic
behavior of the Bessel functions
\begin{equation}\label{Kinf}
    K_\nu(z)\sim\sqrt{\frac{\pi}{2z}}\,\,e^{-z},\qquad\quad
    |z|\to+\infty
\end{equation}
we immediately see that the asymptotic behavior of the \pdf\
\eref{VG1} is $(\alpha|x|)^\lambda e^{-\alpha|x|}$, and hence the
momenta always exist for every $\lambda\in\mathbb{R}$. Of course
this corresponds to the fact that the \chf\ \eref{VG2} is always
derivable in $u=0$. Since our laws are centered and symmetric the
odd momenta vanish; as for the even momenta we have by direct
calculation
\begin{equation*}
    m_{VG}(2k)=\frac{2^k(2k-1)!!}{\alpha^{2k}}\,\frac{\Gamma(\lambda+k)}{\Gamma(\lambda)}\,,\qquad\quad
    k=0,1,2,\ldots
\end{equation*}
so that the expectation is always zero, and the variance
\begin{equation*}
    \sigma^2_{VG}=\frac{2\lambda}{\alpha^2}
\end{equation*}
Then it is easy to see that for a given $\sigma>0$ the laws
$\VG(\lambda,\sqrt{2\lambda}/\sigma)$ have all the same variance
$\sigma^2$ for every value of $\lambda$, and that for
$\lambda\to+\infty$ they converge in distribution to the normal law
$\mathcal{N}(0,\sigma)$. From the \chf\ \eref{VG2} we immediately
see that the VG distributions are \id\ but not stable. It is easy to
see, however, that, the sub--families $\VG(\lambda,\alpha)$ with a
fixed value of $\alpha$ are closed under convolution: in fact the
sum of two independent \rv's respectively with laws
$\VG(\lambda_1,\alpha)$ and $\VG(\lambda_2,\alpha)$ is a \rv\ with
law $\VG(\lambda_1+\lambda_2\,,\alpha)$, as can easily be seen from
\eref{VG2}. This of course does not amount to stability since laws
$\VG(\lambda,\alpha)$ and $\VG(\lambda',\alpha)$ with
$\lambda\neq\lambda'$ are not of the same type. For the sake of
simplicity in the following we will take $\alpha=1$ and we will use
the shorthand notation $\VG(\lambda)=\VG(\lambda,1)$.

\subsection{The Student distributions}

The class of the centered, symmetric Student laws (see Heyde and
Leonenko 2005 for a recent review) can be considered as conjugate to
that of the centered, symmetric VG laws in the sense that here the
roles of the \pdf\ and \chf\ are interchanged. They are the limit
for $\alpha\to0^+$ of the $\GH(\lambda,\alpha,\delta)$ laws with
$\lambda<0$. By taking the new parameter $\nu=-2\lambda>0$, and
recalling that $K_\nu(z)=K_{-\nu}(z)$, the \pdf\ and \chf\ of the
$\GH(\lambda,\alpha,\delta)$ laws become
\begin{eqnarray*}
  f_{GH}(x) &=& \frac{\alpha}{\sqrt{2\pi}}\,\frac{(\delta\alpha)^\frac{\nu}{2}}{ K_\frac{\nu}{2}(\delta\alpha)}\,
                    \frac{K_{\frac{\nu+1}{2}}(\alpha\sqrt{\delta^2+x^2})}{(\alpha\sqrt{\delta^2+x^2})^{\frac{\nu+1}{2}}} \\
 \varphi_{GH}(u) &=& \frac{(\delta\sqrt{\alpha^2+u^2})^\frac{\nu}{2}K_\frac{\nu}{2}(\delta\sqrt{\alpha^2+u^2})}
                            {(\delta\alpha)^\frac{\nu}{2}K_\frac{\nu}{2}(\delta\alpha)}
\end{eqnarray*}
so that from~\Eref{Kzero} in the limit for $\alpha\to0^+$ we get the
\pdf\ and \chf\ of the centered, symmetric Student laws
$\ST(\nu,\delta)$
\begin{eqnarray}
  f_{ST}(x) &=& \frac{1}{\delta\, B\!\left(\frac{1}{2},\frac{\nu}{2}\right)}
                    \left(\frac{\delta^2}{\delta^2+x^2}\right)^\frac{\nu+1}{2}\label{ST1} \\
 \varphi_{ST}(u) &=& 2\,\frac{(\delta|u|)^\frac{\nu}{2}K_\frac{\nu}{2}(\delta|u|)}
                         {2^\frac{\nu}{2}\Gamma\left(\frac{\nu}{2}\right)}\label{ST2}
\end{eqnarray}
where $\nu>0$, $\delta>0$ and $B(z,w)$ is the Beta function
(Abramowitz and Stegun 1968). Here $\delta$ is the scale parameter,
while $\nu$ classifies the different law types. It is also easy to
see that for $|x|\to+\infty$ the Student \pdf\ goes to zero as
$|x|^{-\nu-1}$, so that for a given $\nu$ the moments $m_{ST}(n)$
exist only if $n<\nu$. When they exist, the odd momenta are zero for
symmetry, while the even momenta are
\begin{equation*}
    m_{ST}(2k)=\delta^{2k}\frac{B\left(\frac{1}{2}+k\,,\frac{\nu}{2}-k\right)}{B\left(\frac{1}{2},\frac{\nu}{2}\right)},
               \qquad\quad k=0, 1, 2, \ldots,\quad 2k<\nu
\end{equation*}
In particular the expectation exists (and vanishes) for $\nu>1$,
while the variance exists finite for $\nu>2$ and its value is
\begin{equation}\label{STvar}
    \sigma_{ST}^2=\frac{\delta^2}{\nu-2}
\end{equation}
As a consequence, for $\nu>2$ and for a given $\sigma>0$, the laws
$\ST(\nu,\sigma\sqrt{\nu-2})$ have all the same variance $\sigma^2$,
and it is easy to show that for $\nu\to+\infty$ they converge in
distribution to the normal law $\mathcal{N}(0,\sigma)$. It can be
proved that the Student distributions are \id\ (this is not trivial
at all; see Grosswald 1976a and 1976b, Ismail 1977, Bondesson 1979,
Pitman and Yor 1981, Bondesson 1992), but that they are not stable,
with one notable exception: the $\nu=1$ case, that of the Cauchy
laws $\ST(1,\delta)=\mathcal{C}(\delta)$ which constitute one of the
better known classes of stable laws with \pdf\ and \chf\
\begin{equation*}
    f(x)=\frac{1}{\delta\,\pi}\,\frac{\delta^2}{\delta^2+x^2}\, ,\qquad\quad
    \varphi(u)=e^{-\delta|u|}
\end{equation*}
Besides this case -- and at variance with the VG -- the Student laws
are not even closed under convolution: this makes the study of the
time evolution of a Student process a more complicated and
interesting business which constitutes a relevant part of this
paper. For the sake of simplicity in the following we will take
$\delta=1$ and we will use the shorthand notation
$\ST(\nu)=\ST(\nu,1)$.

\section{The VG process}\label{VGproc}

L\'evy processes produced by means of VG distributions are simple
enough because of their closure under convolution. In fact it is
easy to see from~\eref{VG2} that (taking $\alpha=1$ and $T=1$ to
simplify the notations) for a $\VG(\lambda)$ law the transition
\chf\ of the process (with initial time $s=0$ and position $y=0$) is
\begin{equation}\label{VGchf}
    \Phi(u,t|\lambda)=[\varphi_{VG}(u)]^t
    =\left(\frac{1}{1+u^2}\right)^{\lambda t}
\end{equation}
so that the law of the increment in $[0,t]$ always is a VG law with
the parameter evolving in time; namely, at every $t$, we have
$X(t)\sim\VG(\lambda t)$, and hence the corresponding \pdf\ is
explicitly known at every time and is
\begin{equation}\label{VGpdf}
    p(x,t|\lambda)=\frac{2}{2^{\lambda t}\Gamma(\lambda
    t)\sqrt{2\pi}}\,|x|^{\lambda t-\frac{1}{2}}K_{\lambda t-\frac{1}{2}}(|x|)
\end{equation}
Apparently -- as in the Poisson case -- the laws of the process
belong to the VG family all along the evolution, but this does not
mean that the process is stable since the laws of the VG family are
not of the same type. In fact, with increasing values of $t$, the
distributions of a VG process go throughout all the gamut of the VG
family: what changes with $\lambda$ is just the instant when the
distribution is simply a bilateral exponential. As remarked in the
Section~\ref{ssectionVG} the \pdf~\eref{VGpdf} has an elementary
form only for $t=\frac{1}{\lambda}, \frac{2}{\lambda},\ldots$ but a
great deal of information is available also in the general, non
elementary form. In particular from~\eref{Kzero} and~\eref{Kinf} we
can study the behavior of the \pdf\ both near the origin and in the
asymptotic region. For small $x$ we find
\begin{equation*}
    p(x,t|\lambda)\sim\left\{
                  \begin{array}{ll}
                    |x|^{2\lambda t-1}, &\quad \hbox{for $0<t<\frac{1}{2\lambda}\,$,}  \vspace{5pt}\\
                    -\log |x|, &\quad \hbox{for $t=\frac{1}{2\lambda}\,$,} \vspace{5pt}\\
                    \frac{1}{2\pi}\,\frac{\Gamma\left(\lambda t-\frac{1}{2}\right)}{\Gamma(\lambda t)},
                                  &\quad \hbox{for $\frac{1}{2\lambda}<t\,$,}
                  \end{array}
                \right.
 \qquad\quad    x\to0
\end{equation*}
namely near the origin the \pdf\ has an integrable singularity for
$0<t\leq \frac{1}{2\lambda}$, and thereafter it takes finite values
for $t>\frac{1}{2\lambda}$. As for the asymptotic behavior we have
\begin{equation*}
    p(x,t|\lambda)\sim|x|^{\lambda t-1}\,e^{-|x|},\qquad\quad
    |x|\to+\infty
\end{equation*}
namely it is a negative exponential times a power. It is apparent
then that this asymptotic behavior changes with time since the power
depends on $t$; it is however always dominated by the exponential so
that all the moments exist at every time. From~\Eref{VGpdf} we can
also explicitly calculate the characteristic triplet:
\begin{equation}\label{WVG}
  A=0\,,\qquad B=0\,,\qquad W(z) = \lambda\, \frac{e^{-|z|}}{|z|}
\end{equation}
so that the dimensionless \PIDE\ for the VG process takes the form
\begin{equation*}
    \partial_t p(x,t) =
  \lim_{\epsilon\to0^+}\int_{|z|\geq\epsilon}
                \lambda \,e^{-|z|}\,\frac{p(x+z,t)-p(x,t)}{|z|}\,\,dz
\end{equation*}
A validation of~\eref{WVG} comes then from the L\'evy--Khinchin
formula (Lo\`eve 1987) which here reads
\begin{equation}\label{LKformula}
    \log\varphi(u)=\lim_{\epsilon\to0^+}\int_{|x|\geq\epsilon}
                  \left(e^{iux}-1-\frac{iux}{1+x^2}\right)W(x)\,dx
\end{equation}
and which from~\eref{VG2} and~\eref{WVG} easily reduces itself to
\begin{equation*}
    -\log(1+u^2)=2\int_0^{+\infty}(\cos ux-1)\,\frac{e^{-x}}{x}\,dx
\end{equation*}
a relation which is immediately verified by direct calculation.

\section{The Student process}\label{STproc}

Despite their apparent symmetry and analogy with the VG family, the
processes produced by the Student laws are not so straightforward to
analyze (for recent results about the Student process see Heyde and
Leonenko 2005). The problem is that the Student family
$\ST(\nu,\delta)$ is not even closed under convolution, so that it
is not easy to figure out the general behavior of a Student process
with arbitrary $\nu$. As a consequence we will limit ourselves here
to study the particular case of the $\nu=3$ process whose features
can be fairly understood: this will also give us an insight on the
possible general behavior of these L\'evy processes. It is important
to remark, moreover, that this particular Student process with
$\nu=3$ is the present candidate to describe the increments in the
velocity process for particles in an accelerator beam (Vivoli \etal
2006), and hence its analysis has not a purely academic interest.
Let us introduce now the following notation for the
$\ST(\nu,\delta)$ laws and the corresponding processes: for $\nu>0$
and $\delta>0$
\begin{eqnarray}
  f(x|\,\nu,\delta) &=& f_{ST}(x)=
                  \frac{1}{\delta\, B\!\left(\frac{1}{2},\frac{\nu}{2}\right)}
                    \left(\frac{\delta^2}{\delta^2+x^2}\right)^\frac{\nu+1}{2}\label{STpdf}\\
  f(x|\,\nu) &=& f(x|\,\nu,1)=\frac{1}{B\left(\frac{1}{2},\frac{\nu}{2}\right)}
                                             \left(\frac{1}{1+x^2}\right)^{\frac{\nu+1}{2}}
                 \label{STpdfred}
\end{eqnarray}
so that $f(x|\,\nu)$ from now on will be the \pdf\ of the Student
law $\ST(\nu)=\ST(\nu,1)$ for $\delta=1$. In the same way we can
introduce the reduced form of the \chf\
\begin{eqnarray*}
  \varphi(u|\,\nu,\delta) &=& \varphi_{ST}(u) =
                 2\,\frac{|\delta u|^{\frac{\nu}{2}}\,K_{\frac{\nu}{2}}(|\delta u|)}
                               {2^{\frac{\nu}{2}}\,\Gamma\left(\frac{\nu}{2}\right)} \\
  \varphi(u|\,\nu) &=& \varphi(u|\,\nu,1) = 2\,\frac{|u|^{\frac{\nu}{2}}\,K_{\frac{\nu}{2}}(|u|)}
                               {2^{\frac{\nu}{2}}\,\Gamma\left(\frac{\nu}{2}\right)}
\end{eqnarray*}
Then, by taking $T=1$, the transition \chf\ of the Student process
for the law $\ST(\nu)$ (with initial time $s=0$ and position $y=0$)
is explicitly known and is
\begin{equation*}
  \Phi(u,t|\,\nu) = [\varphi(u|\,\nu)]^t
\end{equation*}
and the corresponding transition \pdf\  is
\begin{equation*}
  p(x,t|\,\nu) =
        \frac{1}{2\pi}\int_{-\infty}^{+\infty}e^{-iux}\Phi(u,t|\,\nu)\,du
        =
        \frac{1}{2\pi}\int_{-\infty}^{+\infty}e^{-iux}[\varphi(u|\,\nu)]^t\,du
\end{equation*}
If we denote as $\ST(\nu,\delta)$--process the Student process such
that its law at $t=T$ is exactly $\ST(\nu,\delta)$ then
$p(x,t|\,\nu)$ will be the \pdf\ of a $\ST(\nu)$--process. In the
following we will perform our calculations on the reduced,
dimensionless quantities only: we can always revert to the
dimensional variables by means of simple transformations. It is easy
to realize from the form of $\Phi(u,t|\,\nu)$ that for $t\to0^+$ the
process approaches a law degenerate in $x=0$, and that along the
evolution of a Student process the marginal $p(x,t|\,\nu)$ no longer
are simple Student \pdf's: after all we know that the Student family
is neither stable, nor closed under convolution. The main problem is
then to find an explicit form for the transition \pdf\ which by
symmetry can be explicitly written as
\begin{equation}\label{processpdf}
    p(x,t|\,\nu) =\frac{1}{\pi}\int_0^{+\infty}\cos(ux)
                 \left[2\,\frac{|u|^{\frac{\nu}{2}}\,K_{\frac{\nu}{2}}(|u|)}
                               {2^{\frac{\nu}{2}}\,\Gamma\left(\frac{\nu}{2}\right)}\right]^t\,du
\end{equation}

\subsection{The Student processes of odd integer index: the $\nu=3$ case}

Since the integration in~\eref{processpdf} can not be performed in
general we will limit ourselves to particular cases. To do that let
us remark that the Student \chf's have an elementary form for odd
integer values of the parameter $\nu$. In fact from \Eref{Ksemint}
we have for $\nu=2n+1$ with $n=0,1,\ldots$ and with $\ell=n-j$
\begin{eqnarray*}
    f(x|\,2n+1)&=&\frac{\Gamma(n+1)}{\sqrt{\pi}\,\Gamma\left(n+\frac{1}{2}\right)}\left(\frac{1}{1+x^2}\right)^{n+1}
    =\frac{(2n)!!}{\pi(2n-1)!!}\left(\frac{1}{1+x^2}\right)^{n+1}
    \\
    \varphi(u|\,2n+1) &=& 2\,\frac{|u|^{n+\frac{1}{2}}\,K_{n+\frac{1}{2}}(|u|)}
                               {2^{n+\frac{1}{2}}\,\Gamma\left(n+\frac{1}{2}\right)}
           = e^{-|u|}\sum_{\ell=0}^n\frac{n!}{(2n)!}\,\frac{(2n-\ell)!}{(n-\ell)!}\,\frac{(2|u|)^\ell}{\ell!}
\end{eqnarray*}
so that the \chf\ is just an exponential times a polynomial in $|u|$
(see~\Tref{STsemint} for a few explicit examples).
\begin{table}
\caption{\label{STsemint}Examples of odd integer order ($\nu=2n+1$),
dimensionless and reduced Student laws.}
\begin{indented}
\item[]\begin{tabular}{@{}llll} \br
$\nu$ & $n$ & $f(x|2n+1)$ & $\varphi(u|2n+1)$ \\
\mr
 1 & 0 & $\frac{1}{\pi}\,(1+x^2)^{-1}$ & $e^{-|u|}$ \vspace{2pt}\\
 3 & 1 & $\frac{2}{\pi}\,(1+x^2)^{-2}$ & $e^{-|u|}(1+|u|)$ \vspace{2pt}\\
 5 & 2 & $\frac{8}{3\pi}\,(1+x^2)^{-3}$ & $e^{-|u|}(1+|u|+\frac{1}{3}|u|^2)$ \vspace{2pt}\\
 7 & 3 & $\frac{16}{5\pi}\,(1+x^2)^{-4}$ & $e^{-|u|}(1+|u|+\frac{2}{5}|u|^2+\frac{1}{15}|u|^3)$ \\
\br
\end{tabular}
\end{indented}
\end{table}
The first case $n=0$, $\nu=1$ is just the stable, reduced Cauchy law
$\mathcal{C}(1)$ which produces the well known Cauchy process. We
can then look at the explicit time evolution of the first non stable
case by taking the $n=1$, $\nu=3$ law, namely the $\ST(3)$--process
with \pdf\
\begin{eqnarray*}
  p(x,t|\,3) &=&\frac{1}{\pi}\int_0^{+\infty}\cos(ux)
                 e^{-tu}(1+u)^t\,du\\
                 &=&\Re\left\{\frac{1}{\pi}\int_0^{+\infty}
                 e^{-(t+ix)u}(1+u)^t\,du\right\}
\end{eqnarray*}
By taking then
\begin{eqnarray*}
  Q(a,z) &=& \frac{1}{\pi}\int_0^{+\infty}e^{-zu}(1+u)^{a-1}du
               =\frac{1}{\pi}\,\frac{e^z}{z^a}\,\Gamma(a,z) \\
  \Gamma(a,z) &=&
  \int_z^{+\infty}e^{-w}w^{a-1}dw\,,\qquad\qquad\qquad
  \Gamma(a,0)=\Gamma(a)
\end{eqnarray*}
where $\Gamma(a,z)$ is the incomplete Gamma function (Abramowitz and
Stegun 1968), we can also write
\begin{equation}\label{pdf3}
    p(x,t|\,3)=\Re\left\{Q(t+1,t+ix)\right\}
    =\Re\left\{\frac{e^{t+ix}\,\Gamma(t+1,t+ix)}{\pi(t+ix)^{t+1}}\right\}
\end{equation}
This new closed form~\eref{pdf3} of the increment laws of the
Student process with $\nu=3$ is now explicitly given for every time
$t>0$: in the following sections we will try to analyze its
properties.

\subsection{Asymptotic behavior of the $\ST(3)$--process}

Since the Student laws are not closed under convolution we know that
$p(x,t|\,3)$ coincides with a Student law only for $t=1$. A first
question is then to check if, that notwithstanding, some important
property of the $t=1$ distribution is preserved along the evolution.
In fact we will see in the following that for an arbitrary fixed,
finite $t>0$ the asymptotic behavior of $p(x,t|\,3)$ for large $x$
is always infinitesimal at the same order $|x|^{-4}$ of the original
$\ST(3)$
\begin{prop}\label{prop01}
If $p(x,t|\,3)$ is the \pdf~\eref{pdf3} of a Student
$\ST(3)$--process, then
\begin{equation*}
    p(x,t|\,3)=\frac{2t}{\pi x^4}+o\left(|x|^{-4}\right),\qquad\quad
    |x|\to+\infty
\end{equation*}
for every given $t>0$.
\end{prop}

\noindent \textbf{Proof}: Let us remember first of all that by
repeated integration by parts of the incomplete Gamma function we
get the following recurrence formula: for a given $a>0$ and
$n=1,2,\ldots$
\begin{eqnarray*}
  Q(a,z) &=& \frac{1}{\pi}\,\frac{e^z}{z^a}\,\Gamma(a,z)
              = \frac{1}{\pi}\sum_{k=0}^{n-1}\frac{\Gamma(a)}{\Gamma(a-k)}\,\frac{1}{z^{k+1}}
                   +R_n(a,z) \\
  R_n(a,z)    &=& \frac{1}{\pi}\,\frac{e^z}{z^a}\,\frac{\Gamma(a)}{\Gamma(a-n)}\,\Gamma(a-n,z)
\end{eqnarray*}
where, from a classical result about this asymptotical expansion
(Gradshteyn and Ryzhik 1980), the remainder $R_n(a,z)$ is an
infinitesimal of order greater than $n$
\begin{equation*}
    |R_n(a,z)| = O\left(|z|^{-n-1}\right),\qquad\quad |z|\to+\infty
\end{equation*}
Then, for $a=t+1$ and $z=t+ix$ with  an arbitrary but fixed $t>0$,
we will have in the limit $|x|\to+\infty$
\begin{equation*}
    |\Re\{R_n(t+1,t+ix)\}|\leq|R_n(t+1,t+ix)|=O\left(|x|^{-n-1}\right)
\end{equation*}
Now take $n=4$: from the previous expansion and~\Eref{pdf3} we have
for $|x|\to+\infty$
\begin{eqnarray*}
  p(x,t|\,3) &=& \Re\{Q(t+1,t+ix)\} \\
                  &=& \frac{1}{\pi}\sum_{k=0}^3\frac{\Gamma(t+1)}{\Gamma(t-k+1)}\,
                            \frac{\Re\left\{(t-ix)^{k+1}\right\}}{(t^2+x^2)^{k+1}}
                                           +o\left(|x|^{-4}\right)
\end{eqnarray*}
while from a direct calculation of the real parts we will find that
the higher powers exactly cancel away from the numerator so that the
leading asymptotic term for $|x|\to+\infty$ is of the order
$|x|^{-4}$; more precisely we have
\begin{eqnarray*}
\fl
 \sum_{k=0}^3\frac{\Gamma( t+1)}{\Gamma( t-k+1)}
                  \,\frac{\Re\left\{( t-i x)^{k+1}\right\}}{( t^2+ x^2)^{k+1}}
 &=&\frac{2 t x^4-4 t^3( t^2-5 t+3) x^2+2 t^5(2 t^2-2 t+1)}{( t^2+ x^2)^4}
 \\
\fl &=& \frac{2 t}{ x^4}+o\left(| x|^{-4}\right)
\end{eqnarray*}
giving finally the statement in our Proposition.
 \hfill$\Box$

\vspace{5pt}\noindent It must be remarked that the previous result
is true for an arbitrary finite, fixed time $t$. For diverging $t$,
however, the reduced law of the process approaches a gaussian: let
$X(t)$ be our $\ST(3)$--process with \pdf\ $p(x,t|\,3)$; then we
know that
\begin{equation*}
    \mathbf{E}[X(t)]=0\,,\qquad\quad \mathbf{Var}[X(t)]=\mathbf{E}\left[X^2(t)\right]=t
\end{equation*}
so that $t^{-1/2}X(t)$ is a centered, reduced \rv\ for every $t$. A
simple look at the \chf's will then shows that in distribution we
have
\begin{equation*}
    \frac{X(t)}{\sqrt{t}}\tod\mathcal{N}(0,1)\,,\qquad t\to+\infty
\end{equation*}
since for large values of $t$ and arbitrary fixed $u$
\begin{equation*}
  [\varphi(u/\sqrt{t}\,|\,3)]^t = \left[e^{-|u|/\sqrt{t}}
                  \left(1+\frac{|u|}{\sqrt{t}}\right)\right]^t \longrightarrow e^{-u^2/2}
                  \,,\qquad t\to+\infty
\end{equation*}

\subsection{The $\ST(3)$--process distribution at integer times $t=n$}

To understand the time evolution of $p(x,t|\,3)$ we can analyze the
forme of this \pdf\ for integral values of the time $t=n=1,2,\ldots$
since in this case the distributions have explicit elementary
expressions. Of course $p(x,n|\,3)$ is nothing else than the
distribution of the sum of $n$ independent $\ST(3)$ \rv's, so that
the following proposition can also be seen as a new result about the
$n$-th convolution of the law $\ST(3)$.
\begin{prop}\label{prop02}
For $n=1,2,\ldots$ we have (within the notations of the present
Section)
\begin{eqnarray*}
  p(x,n|\,3) &=& \sum_{k=0}^nf(x|2k+1,n)\,q_n(k|\,3) \\
  q_n(k|\,3) &=&
  \frac{(-1)^k}{2k+1}\sum_{j=0}^{2k+1}\binom{n}{j}\binom{2k+1}{j}\binom{j}{k}(j+1)!\left(\frac{-1}{2n}\right)^j
\end{eqnarray*}
where $q_n(k|\,3)$ is a discrete probability distribution taking
(strictly) positive values only for $k=1,2,\ldots,n$ (in particular
$q_n(0|\,3)=0$ for every $n$) and such that
\begin{equation*}
    \sum_{k=1}^n\frac{q_n(k|\,3)}{2k-1}=\frac{1}{n}
\end{equation*}
\end{prop}

\noindent \textbf{Proof}:
see~\ref{appendix1}\hfill$\Box$\vspace{5pt}

\begin{figure}
\begin{center}
\includegraphics*[width=13.0cm]{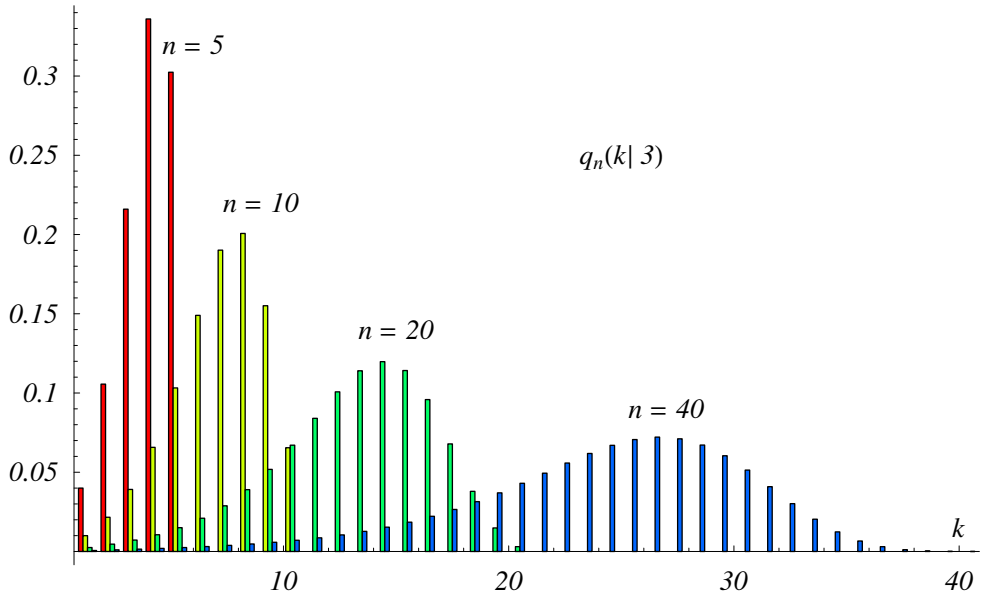}
\end{center}
\caption{Mixture weights of the integer time ($t=n$) components for
a Student process with $\nu=3$.}\label{fig01}
\end{figure}

\noindent The meaning of the Proposition~\ref{prop02} is then that
(at least) at integral times $t=n=1,2,\ldots$ the marginal
one--dimensional \pdf\ $p(x,n|\,3)$ of the $\ST(3)$--Student process
is a mixture (convex combination) of Student \pdf's~\eref{STpdf}
$f(x|\nu,\delta)$ with
\begin{itemize}
    \item odd integer orders $\nu=2k+1$ with $k=0,1,\ldots$,
    \item integer scaling factors $\delta=n$,
    \item relative weights $q_n(k|\,3)$ such that $q_n(0|\,3)=0$, so that
    no Student distribution of order smaller than $\nu=3$ appears in
    the mixture.
\end{itemize}
In other words they are mixtures of $\ST(2k+1,n)$ laws. The
distributions $q_n(k|\,3)$ are a new kind of discrete probability
laws whose bar diagrams at different times $t=n$ are displayed in
Figure~\ref{fig01}.
They show how the weight of the higher order Student distributions
grows with the time, but also that, this notwithstanding, the lowest
order ($\nu=3$) distribution is always present -- albeit with
dwindling importance -- with a non zero weight. We see at once that
this new result is coherent with Proposition~\ref{prop01} and
explicitly shows how the asymptotic behavior is kept $|x|^{-4}$ all
along the time evolution: in fact in the mixture representing
$p(x,n|\,3)$ the lowest order Student distribution always is --
albeit with dwindling weight -- that with $\nu=3$ which
asymptotically behaves as $|x|^{-4}$; all the other components in
the mixture are instead faster infinitesimals. The importance of the
higher orders, however, grows with the time. This is exactly the
behavior recently observed in complex dynamical systems used to
simulate the behavior of intense beams of charged particles in
accelerators (Vivoli et al 2006). Due to their mutual interactions
these particles follow irregular paths, and a statistical analysis
shows that the distribution of the increments follows an almost
gaussian distribution in its central part, and a Student
$\ST(3,\delta)$ distribution on the tails with a $|x|^{-4}$ decay
rate. This suggests that the beam particles follow a $\ST(3,\delta)$
L\'evy process which is observed at a time scale ($\Delta t$) large
when compared to some characteristic time $T$ of the process, but
finite and fixed so that the increment distribution shows two
different regimes (gaussian and $|x|^{-4}$) in the two regions.

The results presented in Proposition~\ref{prop01} and~\ref{prop02}
that the \pdf\ of a L\'evy--Student process is a suitable finite
mixture of other Student \pdf's of different types has been proved
here only in the particular conditions chosen for our demonstration.
It suggests however a possible generalization: it is fair in fact to
put forward the conjecture that every L\'evy--Student process at
every time will have a marginal one dimensional \pdf\ which is a
mixture of other Student \pdf's, but not necessarily (as in our
particular case) of a finite number of odd integer indices Student
\pdf's. In other words, by keeping always the same notation, the
\pdf\ $p(x,t|\nu_0)$ could be a (possibly continuous) mixture of
Student \pdf's $f(x|\,\nu,\delta)$ through a (possibly continuous)
distribution $q_t(\nu|\,\nu_0)$. Finally, in order to preserve the
result of Proposition~\ref{prop01}, we could also conjecture that
$q_t(\nu|\,\nu_0)$ gives probability zero in the mixture to every
Student law with $\nu<\nu_0$. If some form of this conjecture shows
up to be true this would determine some new family of randomized
Student distributions which is closed under convolution.

\subsection{The L\'evy triplet for a $\ST(3)$--process}\label{STproctriplet}

\begin{figure}
\begin{center}
\includegraphics*[width=9.0cm]{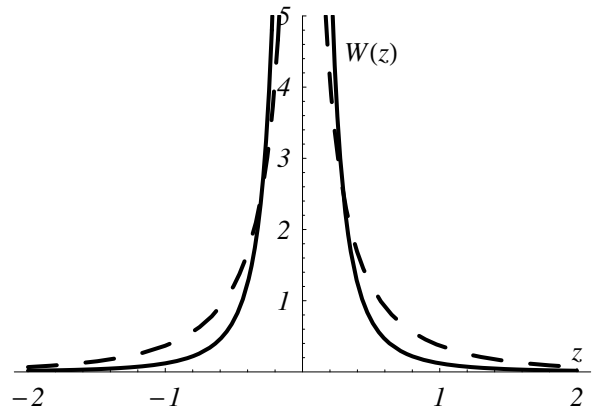}
\end{center}
\caption{Plot of the (reduced and dimensionless) L\'evy densities
for a Student with $\nu=3$ (solid line) and for a VG (dashed line)
process.}\label{fig02}
\end{figure}

We will finally calculate the elements of the L\'evy triplet for a
$\ST(3)$--process from the formulas~\eref{Wchf},~\eref{Achf}
and~\eref{Bchf}. First of all, due to the $\ST(3)$ law symmetry, we
already know that $A=0$; then we must recall that the \chf\ of the
$\ST(3)$ law is
\begin{equation}\label{chf3}
    \varphi(u|\,3)=e^{-|u|}(1+|u|)
\end{equation}
so that by a direct calculation we get an explicit expression of the
L\'evy triplet with $T=1$ (for details on the derivation
see~\ref{appendix2})
\begin{equation}
  A=0\,,\quad B=0\,,\qquad
  W(z)=\frac{1-|z|\left(\,\sin|z|\,\ci|z|-\cos|z|\,\si|z|\,\right)}{\pi
  z^2}\label{W3red}
\end{equation}
where the sine and the cosine integral functions are (Gradshteyn and
Ryzhik 1980)
\begin{equation*}
    \si x= -\int_x^{+\infty}\frac{\sin t}{t}\,dt\,,\qquad\quad\ci x= -\int_x^{+\infty}\frac{\cos t}{t}\,dt
\end{equation*}
A plot of $W(z)$ is shown in Figure~\ref{fig02} where it is also
compared with the analogous density~\eref{WVG} for a VG process.
The behavior of $W(z)$ at the origin and at the infinity is
\begin{equation*}
    W(z)=\left\{
       \begin{array}{ll}
         z^{-2}+o(z^{-2}), & \qquad\hbox{$z\to0^+$;} \\
         2\,z^{-4}+o(z^{-4}), & \qquad\hbox{$z\to+\infty$.}
       \end{array}
     \right.
\end{equation*}
In particular remark that near the origin it has the same behavior
of the L\'evy density for the Cauchy process in~\eref{cauchytriple},
while it asymptotically behaves exactly as the $\ST(3)$
distribution. We could then also conjecture here that the $W(z)$
function of a generic $\ST(\nu,\delta)$--process will always have a
$z^{-2}$ behavior for $z\to0^+$, and a $|z|^{-\nu-1}$ behavior for
$|z|\to+\infty$. From~\eref{W3red} we also see that, always with
$T=1$, the \PIDE~\eref{pide} for a $\ST(3)$--process takes the
particular form
\begin{equation*}
  \partial_t p(x,t) =
  \lim_{\epsilon\to0^+}\int_{|z|\geq\epsilon}W(z)
                \left[p(x+z,t)-p(x,t)\right]\,dz
\end{equation*}
with $W(z)$ given in~\eref{W3red}. Finally, inspection into the
L\'evy--Khinchin formula~\eref{LKformula} for the \chf~\eref{chf3}
immediately gives as a byproduct a previously unknown way to
calculate a non trivial integral:
\begin{equation}
    \frac{2}{\pi}\int_0^{+\infty}
     \frac{\sin z\,\ci z-\cos z\,\si z}{ z}\,(1-\cos u z)\,d z
     =\log(1+|u|)
\end{equation}

\section{Pathwise properties and simulations}\label{pathprop}

Both the classes of processes analyzed in this paper do not have a
Brownian component ($B=0$) in their L\'evy decomposition which is in
fact reduced to its jumping part and has the form (Cont and Tankov
2004, {\O}ksendal and Sulem 2005)
\begin{eqnarray*}
  X(t) &=& \int_{|z|\geq1}zN(t,dz)+\lim_{\epsilon\to0^+}\int_{\epsilon\leq|z|<1}z\tilde{N}(t,dz)\\
  \tilde{N}(t,U) &=& N(t,U)-\frac{t}{T}\,\nu(U)
\end{eqnarray*}
where $U$ is a Borel set $U\subset \mathbb{R}$, $N(t,U)$ is the jump
measure of the process, namely is the number of the (non zero) jumps
of size in $U$ occurring in $[0,t]$, and
$\nu(U)=\mathbf{E}\left[N(1,U)\right]$ is the L\'evy measure of the
process. In fact $N(t,U)$ is a Poisson process of intensity $\nu(U)$
and $\tilde{N}(t,U)$ is the corresponding compensated Poisson
process. The function $W(x)$ introduced in the previous sections of
this paper plays the role of a density for the L\'evy measure in the
sense that $\nu(dx)=TW(x)\,dx$, so that we have all the elements to
characterize the L\'evy decompositions  of our processes. In
particular, due to the nature of the singularities of the $W(x)$
functions in $x=0$, it is possible to see that both the VG and the
$\ST(3,\delta)$ processes (as well as the Cauchy process) have
infinite activity, namely that $\nu(\mathbb{R})=+\infty$. In that
event we know (Cont and Tankov 2004) that the set of jump times of
every trajectory is countably infinite and dense in $[0,+\infty]$.
This property, together with the continuous distributions of the
jump sizes, accounts for the fact that at first sight the
(simulated) samples of both a VG and a $\ST(3,\delta)$ process do
not look very different from that of a Wiener process, in particular
when we compare just the free trajectories of these processes. Then
to better see the respective pathwise characteristics it will be
useful to introduce some L\'evy diffusions, namely the solutions of
other \SDE\ driven by a L\'evy process $X(t)$ (Protter 2004,
Applebaum 2004, {\O}ksendal and Sulem 2005). If $X(t)$ is a pure
jump L\'evy process, let us consider the L\'evy diffusions $Y(t)$
solution of the \SDE
\begin{eqnarray*}
  dY(t) &=& \alpha(t,Y(t))\,dt+dX(t) \\
  dX(t) &=& \int_{|z|\geq1}zN(dt,dz)+\lim_{\epsilon\to0^+}\int_{\epsilon\leq|z|<1}z\tilde{N}(dt,dz)
\end{eqnarray*}
which is nothing else than a deterministic dynamic system
$\dot{y}(t)=\alpha(t,y(t))$ perturbed by a jump noise $X(t)$. The
simplest case is that of a linear force $\alpha(y)=-ky$ giving rise
to non--Gaussian Ornstein--Uhlenbeck (OU) processes (see for example
Barndorff--Nielsen and Shephard 2001, Cont and Tankov 2004)
\begin{equation}\label{OUlevysde}
    dY(t) = -k\,Y(t)\,dt+dX(t)
\end{equation}
The usual, Gaussian OU process, on the other hand, is the solution
of a \SDE\ where the noise $B(t)$ is completely Brownian with no
jump component:
\begin{equation}\label{OUsde}
    dY(t) = -k\,Y(t)\,dt+dB(t)
\end{equation}
We can compare now the samples of OU--type processes driven either
by a Brownian noise, or by a pure jump noise as the VG and the
$\ST(3,\delta)$ processes. To do that we will produce samples of
$5\,000$ steps by using reduced and dimensionless versions of our
distributions that we will take of unit variance. In particular we
will suppose that for time intervals $\Delta t =T$ the laws of the
noise increments are that reproduced in Table~\ref{incrementlaws}.
\begin{table}
\caption{\label{incrementlaws}The unit variance laws and \pdf's of
the increments used in producing the samples of Figure~\ref{fig03}.}
\begin{indented}
\item[]\begin{tabular}{@{}ccc} \br
 (a) & (b) & (c) \\
\mr
 $\mathcal{N}(0,1)$ & $\VG(1,\sqrt{2})$ & $\ST(3,1)$ \vspace{5pt} \\
 $\frac{1}{\sqrt{2\pi}}\,e^{-x^2/2}$ & $\frac{1}{\sqrt{2}}\,e^{-\sqrt{2}\,|x|}$ &
           $\frac{2}{\pi}\,\frac{1}{(1+x^2)^2}$ \\
\br
\end{tabular}
\end{indented}
\end{table}
Of course the choice of $\Delta t=T$ is instrumental because the VG
and the Student laws have distributions of elementary form only for
$\Delta t=nT$ with $n$ integer (and particularly simple for $n=1$),
as we have seen in the previous Sections. It is not so easy, on the
other hand, to produce our pure jump driven trajectories at other
time scales, in particular for time scales which are fractions of
$T$. At first sight we could think to overcome this difficulty by
arbitrarily changing the value of $T$, but we should remember from
our previous discussion (Section~\ref{lproc}) that our pure jump
processes are not scale invariants, so that different values of $T$
produce different processes.
\begin{figure}
\begin{center}
\includegraphics*[width=13.0cm]{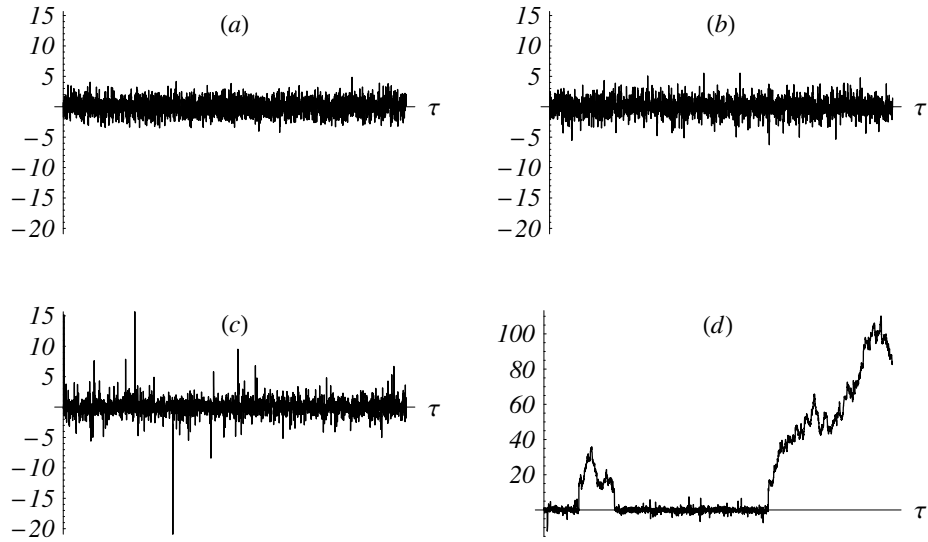}
\end{center}
\caption{Samples of OU--type diffusions ($\tau=t/T$): (a) usual OU
process driven by gaussian Brownian motion; (b) OU--type process
driven by a VG L\'evy noise; (c) OU--type process driven by a
Student L\'evy noise; (d) OU--type process with Student noise and
restoring force of limited range.}\label{fig03}
\end{figure}
Examples of simulated samples of these processes are produced by
discretizing our \SDE\ and are shown in Figure~\ref{fig03} as
functions of the dimensionless time $\tau=t/T$. The parts (a), (b)
and (c) show trajectories produced by our three different \SDE's:
while (a) is a typical sample of an OU process solution of the
\SDE~\eref{OUsde} driven by a normal Brownian motion, the parts (b)
and (c) display typical trajectories produced by the
\SDE~\eref{OUlevysde} driven by respectively a VG noise and a
Student noise. The plots are on the same spatial scale and we can
see the jumping nature of the non gaussian noises from the fact
that, while trajectory (a) is rather strictly confined inside the
region determined by the restoring force $-ky$, the trajectory (b),
and above all the trajectory (c) show clearly random spikes going
outside the confining region. These spikes are produced by the jumps
of the driving processes, and the fact that the Student (c) spikes
are larger than that of the VG (b) case depends on the fact that the
VG distribution has exponential tails which -- albeit longer than
the gaussian tails -- are much shorter than the power tails of a
Student distribution (see also the corresponding asymptotic behavior
of the L\'evy densities $W(z)$ displayed in the
Sections~\ref{VGproc} and~\ref{STproctriplet}). The size of the
spikes can also be put in evidence by cutting the restoring force of
the \SDE's to a finite length, namely by considering the solutions
of
\begin{eqnarray*}
  dY(t) &=& \alpha(Y(t))\,dt+dX(t) \\
  \alpha(y)&=&\left\{
                  \begin{array}{ll}
                    -ky, & \hbox{for $|y|\leq q$;} \\
                    0, & \hbox{for $|y|> q$.}
                  \end{array}
                \right.
                \qquad\quad q>0
\end{eqnarray*}
In this case the restoring force acts only when the process lies in
$[-q,q]$, while the process is completely free outside this region.
Hence when the process jumps beyond the boundaries in $y=\pm q$ it
begins to diffuse freely drifting away from the bounding region.
Occasionally, however, it can also be recaptured by the binding
force. All these features are represented in the part (d) of
Figure~\ref{fig03} which displays the trajectory of a Student driven
OU--type process with a limited range of the force. To compare it
with the other two cases we must now look at the different values of
$q$ that make an escape reasonably likely: while to let an OU
gaussian process to escape is necessary to have a rather small value
of $q$, evasions are likely in the VG case for larger, and in the
Student case even for much larger, $q$ values.

\section{Conclusions}

We have studied in this paper a few examples of non stable,
infinitely divisible processes, and in particular we have explicitly
written down their evolution equations and the laws of the
increments which are the germ of the corresponding markovian
evolutions. In particular we focused our attention on the Student
processes and we presented a new explicit form of their transition
functions. Since the Student family of laws is infinitely divisible,
but non closed under convolution the distribution of the
corresponding L\'evy--Student process is a Student distribution only
at the one particular time. Along the evolution, instead, the
process distribution is no longer a simple Student distribution. We
have shown in the previous sections that, this notwithstanding, at
least in the case of a specific type of Student distribution (with
finite variance), and at least in an infinite sequence of
equidistant time instants the process transition law is a mixture of
a finite number of Student laws given by means of a new kind of
discrete probability distribution. This prompts the conjecture that
in fact while the Student family is not closed under convolution,
some family of mixtures of Student distributions can possibly be
closed. On the other hand, while it is easy to show that for large
values of time the reduced increment law tends to be normal (as it
should be since we are dealing with finite variance distributions),
we have also emphasized that for a finite (albeit large) time the
asymptotic behavior always is the same as that of the Student
distribution at the unit time. This behavior has been put in
evidence by Vivoli \etal 2006 in their model for halo in particle
beams, and we have put forward the conjecture that this could also
be a more general behavior of the Student processes. This last
remark is interesting also in connection with a possible
generalization of the stochastic mechanics that we mentioned in the
Section~\ref{intro} and that will be the argument of forthcoming
research.

\ack The author want to thank C.\ Benedetti, F.\ Mainardi, G.\
Turchetti and A.\ Vivoli for useful discussions and suggestions, and
F.\ De Martino, S.\ De Siena, F.\ Illuminati and M.\ Pusterla for
the long collaboration which was -- and is -- instrumental in this
inquiry.

\appendix

\section{Proof of Proposition~\ref{prop02}}\label{appendix1}
For $a=n+1$ and $n=1,2,\ldots$ the incomplete Gamma functions have a
finite elementary expression (Gradshteyn and Ryzhik 1980) so that
\begin{equation*}
    Q(n+1,z)=\frac{e^z}{z^{n+1}}\,\Gamma(n+1,z)=\sum_{j=0}^n\frac{n!}{(n-j)!}\,\frac{1}{z^{j+1}}
\end{equation*}
and hence we get
\begin{eqnarray*}
  p(x,n|\,3) &=& \frac{1}{\pi}\,\Re\{Q(n+1,n+i x)\}
                  =\frac{1}{\pi}\,\sum_{j=0}^n\frac{n!}{(n-j)!}\,\Re\left\{\frac{1}{(n+i x)^{j+1}}\right\} \\
   &=&\frac{1}{\pi}\,\sum_{j=0}^n\frac{n!}{(n-j)!}\,\frac{1}{(n^2+ x^2)^{j+1}}\,
               \Re\left\{\sum_{m=0}^{j+1}\binom{j+1}{m}(-i x)^mn^{j-m+1}\right\}\\
   &=&\frac{1}{\pi}\,\sum_{j=0}^n\frac{n!}{(n-j)!}\,\frac{1}{(n^2+ x^2)^{j+1}}
               \sum_{2\ell=0}^{j+1}\binom{j+1}{2\ell}(-1)^\ell  x^{2\ell}n^{j-2\ell+1}
\end{eqnarray*}
where it is understood that the second sum is extended to all the
integer values of $\ell$ such that $0\leq2\ell\leq j+1$, namely: if
$j$ is even then $\ell=0,1,\ldots,\frac{j}{2}$; if $j$ is odd then
$\ell=0,1,\ldots,\frac{j+1}{2}$. A little manipulation and the use
of~\Eref{STpdf} then give
\begin{eqnarray*}
  p( x,n|\,3) &=& \frac{1}{\pi}\,\sum_{j=0}^n\binom{n}{j}\frac{j!}{n^{j+1}}\left(\frac{n^2}{n^2+ x^2}\right)^{j+1}
                      \sum_{2\ell=0}^{j+1}\binom{j+1}{2\ell}(-1)^\ell\left(\frac{ x^2}{n^2}\right)^\ell \\
  &=&
  \frac{1}{\pi}\,\sum_{j=0}^n\binom{n}{j}\frac{j!}{n^j}\sum_{2\ell=0}^{j+1}\binom{j+1}{2\ell}\!
           \sum_{m=0}^\ell\binom{\ell}{m}\frac{(-1)^m}{n}\left(\frac{n^2}{n^2+ x^2}\right)^{j-m+1}\\
  &=& \frac{1}{\pi}\,\sum_{j=0}^n\binom{n}{j}\frac{j!}{n^j}\sum_{2\ell=0}^{j+1}\binom{j+1}{2\ell}
           \sum_{m=0}^\ell(-1)^m\binom{\ell}{m}\\
           & & \qquad\qquad\qquad\times B\left(\frac{1}{2}\,,\,j-m+\frac{1}{2}\right) f( x|\,2(j-m)+1,\,n)
\end{eqnarray*}
with $f(x|\nu)$ defined in~\eref{STpdf}. Now by exchanging the order
of the last two sums (with the previous conventions about the range
of the indexes $\ell$ and $m$) we have with $k=j-m$
\begin{eqnarray*}
\fl  p( x,n|\,3) &=&
       \frac{1}{\pi}\,\sum_{j=0}^n\binom{n}{j}\frac{j!}{n^j}\sum_{2m=0}^{j+1}(-1)^m
                            f( x|\,2(j-m)+1,\,n) \\
\fl   & & \qquad\qquad\qquad\qquad\quad\times
                 B\left(\frac{1}{2}\,,\,j-m+\frac{1}{2}\right)
                            \sum_{2\ell=2m}^{j+1}\binom{j+1}{2\ell}\binom{\ell}{m}  \\
\fl   &=&
              \frac{1}{\pi}\,\sum_{j=0}^n\binom{n}{j}\frac{j!}{n^j}\sum_{2m=0}^{j+1}(-1)^m
                            f( x|\,2(j-m)+1,\,n) \\
\fl   & & \qquad\qquad\qquad\qquad\quad\times
            B\left(\frac{1}{2}\,,\,j-m+\frac{1}{2}\right)
                            \frac{2^{j-2m}(j+1)(j-m)!}{m!(j-2m+1)!}  \\
\fl &=& \sum_{j=0}^n\binom{n}{j}\frac{1}{(2n)^j}\sum_{2m=0}^{j+1}
                   \frac{(-1)^m(j+1)!(2j-2m)!}{(j-m)!m!(j-2m+1)!}\,f( x|\,2(j-m)+1,\,n) \\
\fl &=& \sum_{j=0}^n\binom{n}{j}\frac{1}{(2n)^j}\sum_{2k\geq
                                             j-1}^{2j}
                 \frac{(-1)^{j-k}(j+1)!(2k)!}{k!(j-k)!(2k-j+1)!}\,f( x|\,2k+1,\,n)
\end{eqnarray*}
where it is understood that the second sum extends over all the $k$
values such that $j-1\leq2k\leq2j$, namely: for odd $j$ we have
$k=\frac{j-1}{2},\ldots,j$, while for even $j$ we have
$k=\frac{j}{2},\ldots,j$. Finally, by exchanging again the sums and
by adopting the convention that a binomial symbol $\binom{a}{b}$
always is zero whenever the limitation $b\leq a$ is not verified, we
have the results of Proposition~\ref{prop02}
\begin{eqnarray*}
  p(x,n|\,3) &=&\sum_{k=0}^nf(x|\,2k+1,n)\,q_n(k|\,3) \\
  q_n(k|\,3)&=&\frac{(-1)^k}{2k+1}\sum_{j=0}^{2k+1}
             \binom{n}{j}\binom{2k+1}{j}\binom{j}{k}(j+1)!\left(\frac{-1}{2n}\right)^j
\end{eqnarray*}
Since the distribution of our Student process is now represented as
a linear combination of the Student $\ST(2k+1,n)$ \pdf's,
$p(x,n|\,3)$  turns out to be a (randomized, Feller 1971) mixture,
and the coefficient $q_n(k|\,3)$ of this combination must satisfy
\begin{equation*}
    q_n(k|\,3)\geq0\,,\qquad\qquad\sum_{k=0}^nq_n(k|\,3)=1
\end{equation*}
with $q_n(0|3)=0$ for every $n$, as can be seen by direct
calculation. Hence we have that $q_n(k|\,3)$ is a discrete
probability distribution taking non--zero values only for
$k=1,2,\ldots,n$. Finally by remembering that our Student process
has zero expectation and variance $t=n$, and taking also into
account the~\Eref{STvar}, we can write
\begin{eqnarray*}
    n&=&\int_{-\infty}^{+\infty}x^2p(x,n|\,3)\,dx
        =\sum_{k=0}^nq_n(k|\,3)\int_{-\infty}^{+\infty}x^2f(x|2k+1,n)\,dx\\
     &=&\sum_{k=0}^nq_n(k|\,3)\frac{n^2}{2k-1}
\end{eqnarray*}
so that we immediately get also the last result in our proposition.

\section{Derivation of equation~\eref{W3red}}\label{appendix2}
From~\eref{Bchf} and~\eref{chf3} we have for a $\ST(3)$--process
that
\begin{eqnarray*}
  B &=& \frac{1}{\pi }\,\lim_{\epsilon\to0^+}\lim_{M\to+\infty}
                    \int_{-M}^M\frac{-u}{1+|u|}\,\frac{u\epsilon\cos u\epsilon-\sin
                    u\epsilon}{u^2}\,du \\
   &=& \frac{2}{\pi }\,\lim_{\epsilon\to0^+}
                    \int_0^{+\infty} \frac{\sin
                    u\epsilon-u\epsilon\cos u\epsilon}{u(1+ u)}\,du\\
   &=& \frac{2}{\pi }\,\lim_{\epsilon\to0^+}
        \left[\frac{\pi}{2}-\left(\ci\epsilon
                        -\epsilon\,\,\si\epsilon\right)\sin\epsilon
        +\left(\epsilon\,\,\ci\epsilon
                                   +\si\epsilon\right)\cos\epsilon\right]=0
\end{eqnarray*}
where the sine and the cosine integral functions are defined in the
text: hence, as for the Cauchy and the VG processes, the Brownian
part is absent also in this Student process. As for the L\'evy
density $W(z)$, from~\eref{Wchf} we get
\begin{eqnarray*}
  W(z) &=&  \frac{1}{2\pi iz}\lim_{M\to+\infty}\int_{-M}^M\frac{-u}{1+|u|}\,e^{-iuz}\,du\\
   &=& \frac{1}{\pi |z|}\lim_{M\to+\infty}\int_0^{M}\frac{u}{1+u}\sin(u|z|)\,du \\
   &=& \frac{1+|z|\left(\cos|z|\,\,\si|z|-\sin|z|\,\,\ci|z|\right)}{\pi z^2}
\end{eqnarray*}
so that for our $\ST(3)$--process we finally have~\eref{W3red}.

\References

\item[] Abramowitz M and Stegun I A 1968 \textit{Handbook of Mathematical Functions}
(New York: Dover Publications)

\item[] Albeverio S, Blanchard P and H\o gh-Krohn R 1983 \textit{Expo. Math.} \textbf{4} 365

\item[] Albeverio S, R\"udiger B and Wu J--L 2001 in \textit{L\'evy processes, Theory and
applications} ed Barndorff--Nielsen O \etal (Boston: Birkh\"auser)
p~187

\item[] Applebaum D 2004 \textit{L\'evy processes and stochastic
calculus} (Cambridge: Cambridge University Press)

\item[] Barndorff--Nielsen O E 2000 Probability densities and L\'evy
densities (MaPhySto, Aarhus, \textit{Preprint} MPSRR/2000-18)

\item[] Barndorff--Nielsen O E and Shephard N 2001 in
\textit{L\'evy processes, Theory and applications} ed
Barndorff--Nielsen O \etal (Boston: Birkh\"auser) p~283

\item[] Barndorff--Nielsen O E, Mikosch T and Resnick S I (eds) 2001 \textit{L\'evy processes, Theory and
applications} (Boston: Birkh\"auser)

\item[] Bondesson L 1979 \textit{Ann. Prob.} \textbf{7} 965

\item[] Bondesson L 1992 \textit{Generalized Gamma Convolutions and Related
Classes of Distributions and Densities} (Lecture Notes in Statistics
vol 76; Berlin: Springer)

\item[] Bouchaud J--P and Georges A 1990 \textit{Phys. Rep.} \textbf{195} 127

\item[] Chechkin A V, Gonchar V Y, Klafter J, Metzler R and
Tanatarov L V 2004 \textit{J.\ Stat.\ Phys.}\ \textbf{115} 1505

\item[] Chechkin A V, Klafter J, Gonchar V Y, Metzler R and
Tanatarov L V 2005 \PR \textit{E} \textbf{67} 010102(R)

\item[] Cont R and Tankov P 2004 \textit{Financial modelling with jump
processes} (Boca Raton: Chapman\&Hall/CRC)

\item[] Cufaro Petroni N, De Martino S, De Siena S and Illuminati F
1999 \JPA \textbf{32} 7489

\item[] Cufaro Petroni N, De Martino S, De Siena S and Illuminati F
2000 \PR \textit{E} \textbf{63} 016501

\item[] Cufaro Petroni N, De Martino S, De Siena S and Illuminati F
2003 \PR \textit{ST Accelerators and Beams} \textbf{6} 034206

\item[] Cufaro Petroni N, De Martino S, De Siena S and Illuminati F
2004  \textit{Int. J. Mod. Phys. B} \textbf{18} 607

\item[] Cufaro Petroni N, De Martino S, De Siena S and Illuminati F
2005 \PR \textit{E} \textbf{72} 066502

\item[] Cufaro Petroni N, De Martino S, De Siena S and Illuminati F
2006 \NIM \textit{A} \textbf{561} 237

\item[] De Angelis G F 1990 \JMP \textbf{31} 1408

\item[] De Angelis G F and Jona--Lasinio G 1982 \JPA \textbf{15} 2053

\item[] Eberlein E 2001 in
\textit{L\'evy processes, Theory and applications} ed
Barndorff--Nielsen O \etal (Boston: Birkh\"auser) p~371

\item[] Eberlein E and Raible S 2000 \textit{European Congress
of Mathematics (Barcelona)} vol II (\textit{Progress in Mathematics}
vol 202) ed Casacuberta C \etal (Basel: Birkh\"auser) p~367

\item[] Feller W 1971 \textit{An introduction to probability theory and its
applications} vol II (New York: Wiley\& Sons)

\item[] Gardiner C W 1997 \textit{Handbook of stochastic methods}
(Berlin: Springer)

\item[] Gnedenko B V and Kolmogorov A N 1968 \textit{Limit distributions for
sums of independent random variables} Reading: Addison--Wesley)

\item[] Gorenflo R and Mainardi F 1998a \textit{Fract.\ Calc.\ Appl.\
An.}\ \textbf{1} 167 (reprinted at http://www.fracalmo.org/)

\item[] Gorenflo R and Mainardi F 1998b \textit{Arch.\ Mech.}\
\textbf{50} 377 (reprinted at http://www.fracalmo.org/)

\item[] Gradshteyn I S and Ryzhik I M 1980 \textit{Table of integrals,
series  and products} (San Diego: Academic Press)

\item[] Grosswald E 1976a \textit{Ann.\ Prob.}\ \textbf{4} 680

\item[] Grosswald E 1976b \textit{Z.\ Wahrsch.}\ \textbf{36} 103

\item[] Guerra F 1981 \textit{Phys. Rep.} \textbf{77} 263

\item[] Guerra F and Morato L M 1983 \PR \textit{D} \textbf{27} 1774

\item[] Heyde C C and Leonenko N N 2005 \textit{Adv. Appl. Prob.} \textbf{37}
342

\item[] Ishikawa Y 1994 \textit{Tohoku Math.\ J.}\ \textbf{46} 443

\item[] Ismail M E H 1977 \textit{Ann.\ Prob.}\ \textbf{5} 582

\item[] Jacob N and Schilling R L 2001 in \textit{L\'evy processes, Theory and
applications} ed Barndorff--Nielsen O \etal (Boston: Birkh\"auser)
p~139

\item[] L\'eandre R 1987 Densit\'e en temps petit d'un processus de
sauts \textit{S\'eminaire de Probabilit\'es XXI} (\textit{Lecture
Notes in Mathematics} vol 1247) ed J Az\'ema, P A Meyer and M Yor
(Berlin, Springer) p~81

\item[] Lo\`eve M 1987 \textit{Probability theory} vol I (New York:
Springer)

\item[] Lo\`eve M 1978 \textit{Probability theory} vol II (New York:
Springer)

\item[] Madan D B, Carr P P and Chang E C 1998 \textit{European Finance
Review} \textbf{2} 79

\item[] Madan D B and Milne F 1991 \textit{Mathematical Finance} \textbf{1(4)} 39

\item[] Madan D B and Seneta E 1987 \textit{Journal of the Royal Statistical
Society} series B \textbf{49(2)} 163

\item[] Madan D B and Seneta E 1990 \textit{Journal of Business} \textbf{63}
511

\item[] Mantegna R and Stanley H E 2001 \textit{An introduction to
econophysics} (Cambridge: Cambridge University Press)

\item[] Metzler R and Klafter J 2000 \textit{Phys. Rep.} \textbf{339} 1

\item[] Morato L 1982 \JMP \textbf{23} 1020

\item[] Nelson E 1967 \textit{Dynamical theories of Brownian motion}
(Princeton: Princeton University Press)

\item[] Nelson E 1985 \textit{Quantum Fluctuations} (Princeton: Princeton University Press)

\item[] {\O}ksendal B and Sulem A 2005 \textit{Applied stochastic
control of jump diffusions} (Berlin: Springer)

\item[] Paul W and Baschnagel J 1999 \textit{Stochastic processes:
from physics to finance} (Berlin: Springer)

\item[] Pitman J and Yor M 1981 in \textit{Stochastic Integrals}
(Lecture Notes in Mathematics vol 851) ed Williams D (Berlin:
Springer) p~285

\item[] Protter P 2005 \textit{Stochastic integration and differential
equations} (Berlin, Springer)

\item[] Raible S 2000 L\'evy processes in finance: theory, numerics
and empirical facts, PhD Thesis (Freiburg University)

\item R\"uschendorf L and Woerner J H C 2002 \textit{Bernoulli}
\textbf{8} 81

\item[] Sato K 1999  \textit{L\'evy processes and infinitely divisible
distributions} (Cambridge, CUP)

\item[] Vivoli A, Benedetti C and Turchetti G 2006 \NIM A \textbf{561}
320

\item[] Woyczy\'nski W A 2001 in \textit{L\'evy processes, Theory and
applications} ed Barndorff--Nielsen O \etal (Boston: Birkh\"auser)
p~241

\endrefs

\end{document}